\documentclass{article}
\usepackage{a4}
\usepackage{amssymb}
\usepackage{mathrsfs}

\input xypic
\input xy
\xyoption{all}

\newtheorem{theorem}{Theorem}[section]
\newtheorem{lemma}[theorem]{Lemma}

\newtheorem{definition}[theorem]{Definition}

\newcommand{\dom}{\mathrm{dom}}

\newcommand{\PERs}{\mathcal P}

\newenvironment{proof}{\noindent {\bf Proof}. \nopagebreak
}{\nopagebreak\hfill$\Box$ \medskip}

\title{A Note on ``Extensional PERs''}
\author{W.P.Stekelenburg}
\newcommand{\N}{\mathbb N}

\begin{document}

\maketitle

\begin{abstract}
The paper ``Extensional PERs'' by P.Freyd, P.Mulry, G.Rosolini and D.Scott (\cite{exp}) identifies a reflective subcategory of the category of PERs, namely the category $\mathcal C$ of pointed complete extensional PERs, which has the interesting property of being \emph{algebraically compact} with respect to \emph{realizable} functors. Unfortunately, the definition of realizable functors used in \cite{exp} is too restrictive, and this is a problem, because specifically that part of the definition that is too restrictive, is a necessary premise to the given algebraic compactness proof. Here, I present two ways to bypass this problem, and thus to complete the proof.
\end{abstract}

The paper ``Extensional PERs'' by P.Freyd, P.Mulry, G.Rosolini and D.Scott (\cite{exp}) identifies a reflective subcategory of the category of PERs, namely the category $\mathcal C$ of pointed complete extensional PERs, which has the interesting property of being \emph{algebraically compact}. Algebraic compactness ensures the existence of solutions to recursive domain equations (see \cite{acc}). In other words: given a functor $F:(\mathcal C^{op})^n\times\mathcal C^m\rightarrow \mathcal C$, there is a \emph{fixobject} $X$ with a \emph{fixmap} $f:F(\vec X)\rightarrow X$, which is an isomorphism. Due to this property, $\mathcal C$ is an interesting candidate for a categorical semantics of programming languages with recursively defined types.
 
There is one restriction though: the functor $F$ has to be \emph{realizable}. The category of PERs and this subcategory of \emph{pointed CEPERs}, are internal in the effective topos. Any internal functor between these categories comes with a realizer for its functorial properties. Hence the name `realizable functor'. Unfortunately, the definition given in the paper seems to be more restrictive: it is not clear to me that all internal functors are realizable according to the definition found in the paper. And this is a problem, because specifically that part of the definition that is too restrictive, is a necessary premise to the given algebraic compactness proof.

In the research for my master thesis I found two ways to bypass this problem. Firstly, weakly complete internal categories, like the category of PERs and the category of pointed CEPERs, already satisfy the weaker property of algebraic completeness. Secondly, any internal functor is isomorphic to some other internal functor that \emph{does} satisfy the more restrictive definition, showing that the original proof can be used without loss of generality.

\section{The Category of PERs: Notation}
In stead of using the notation of \cite{exp}, I will write about PERs with a more usual mathematical symbolism. So:
\begin{definition}\textnormal{
A \emph{PER} is a \emph{partial equivalence relation} on the natural numbers. So a PER $R$ is a subset of $\N^2$ such that:
\begin{itemize}
\item for all $(n,m)\in\N$, $(m,n)\in\N$ (\emph{symmetry})
\item for all $(n,m)$ and $(m,p)\in\N$, $(n,p)\in\N$ (\emph{transitivity})
\end{itemize}
}\end{definition}

Any PER $R$ forms a total equivalence relation on its \emph{domain} $\dom R:= \{n|(n,n)\in R\}$, and the quotients $\dom R/R$ are used to define morphisms between PERs. Given $n\in\dom R$, I use $[n]_R$ to denote the equivalence class containing $n$ in $\dom R/R$:
\begin{definition}\textnormal{
A \emph{morphism of PERs} $f:R\to S$, is a function $f:\dom R/R\to\dom S/S$, which is \emph{tracked} by a partial recursive function. This means that there is a partial recursive function $\phi$ such that for all $n\in\dom R$ $\phi n$ is defined and $f[n]_R = [\phi n]_S$.
}\end{definition}

These objects and morphisms form the category of PERs $\PERs$. This category is Cartesian closed, basically because we can define $S^R$ to be the PER of indices of tracking partial recursive functions. So any $f:R\to S$ can be identified with the set of those natural numbers that are indices of tracking functions of $f$. Therefore I will sometimes use $[n]_{R\to S}$ to refer to the function $R\to S$ that is tracked by the $n$-th partial recursive function.
Finally, I write the application of the $n$-th partial recursive function to some number $m$ as a simple juxtaposition: $nm$.

\section{Realizable and Monotone Functors} 
The proper definition of realizable functors, based on the idea that they are internal functors in the effective topos is:

\begin{definition}\textnormal{
An endofunctor $F$ of the category of PERs is \emph{realizable}, if there is a single partial recursive function $\phi$ that tracks $F:\hom(R,S)\to \hom(FR,FS)$ for all $R$ and $S$. This means that $\phi x$ converges whenever $[x]:R\to S$ for every pair of PERs $R$ and $S$, and that \begin{equation}F([x]_{R\to S}) = [\phi x]_{FR\to FS}\label{relfunc}\end{equation}
We say that $\phi$ tracks $F$ is this case.
}\end{definition}

The definition is similar to the definition found in \cite{exp}. What is left out, is the requirement that for some index $i$ of the identity function on $\N$, $\phi i=i$. Because $F$ preserves identities, and because $i$ tracks the identity function on any PER $R$, we know that $F([i]_{R\to R}) = [\phi i]_{FR\to FR} = [i]_{FR\to FR}$. So $i\in\bigcap_R F([i]_{R\to R})$ \emph{does} hold. This still doesn't guarantee that $\phi i = i$, however.

Let $\psi i = i$ and $\psi x = \phi x$ if $x\neq i$. $\psi$ is a recursive function, and one might wonder if it can take the place of $\phi$, saving the original definition. Obviously, ($\ref{relfunc}$) is satisfied for $x\neq i$. In the case that $S=R$, the same equation holds for $i$. So we're left with the case $x=i$ and $S\neq R$.

Now note that $[i]:R\to S$ iff $R\subset S$. Therefore if $R\subset S$, and if $F([i]_{R\to S}) = [\psi i]_{FR\to FS} = [i]_{FR\to FS}$, then $FR\subset FS$. This means that all functors which are tracked by an $i$ preserving function are monotone mappings of PERs. On the other hand, if a functor is monotone this way, and has a tracking function $\phi$, the function $\psi$ defined above is another tracking function, and this one is $i$ preserving. So the functors defined in \cite{exp} are a special kind of realizable functor:

\begin{definition}\textnormal{
A realizable endofunctor $F$ of the category of PERs is \emph{monotone}, if its object map is monotone with respect to the inclusion ordering on PERs. In other words: if $R\subset S$, then $FR\subset FS$.
}\end{definition}


I have not been able to prove (or refute, by the way) that all realizable functors are monotone, or to find a proof in the literature. Sadly, in \cite{exp} the least fixpoints that monotone functors have, are used in the algebraic compactness proof: for any monotone functor $F$ we have a fixpoint $X:= \bigcap\{R|FR\subset R\}$, where $FX=X$ and $i$ represents the fixmap.





\section{Algebraic Completeness}
I'll start with the following general result:
\begin{lemma}
Given any topos $\mathcal E$, and any complete internal category $\mathcal C$. $\mathcal C$ is algebraically complete, meaning: for any internal endofunctor $F$, there is an initial algebra.
\end{lemma}

\begin{proof}
$\mathcal E$ allows the construction of the category of algebras of any endofunctor $F$ of $\mathcal C$ internally, so both the category of $F$-algebras $F$-alg, and the underlying object functor $U:F\textrm{-alg}\to\mathcal C$ are internal to $\mathcal E$. Now this underlying functor creates limits, and since $\mathcal C$ is complete (relative to $\mathcal E$), $F$-alg must be complete too. Therefore it has an initial object, which is an initial algebra for $F$.
\end{proof}

In this general proof we actually only need the limit over one functor, namely the underlying object functor $U:F\textrm{-alg}\to\mathcal C$. The category of PERs $\PERs$ is a \emph{weakly} complete full internal subcategory of the effective topos. Weak completeness means that although for arbitrary internal categories $\mathcal D$ and internal functors $F:\mathcal D\to\mathcal C$ a limiting cone exists, there is no internal functor $\mathcal C^{\mathcal D}\to\mathcal C$ adjoint to the functor $K:\mathcal C\to\mathcal C^{\mathcal D}$ that maps objects to constant functors. This, however, still suffices to prove a limit exists for the underlying PER functor $U:F\textrm{-alg}\to\PERs$, and that weakly complete categories are also algebraically complete.

When we work out the construction of this limit, which is supposed to be a subobject of the product of all algebras, we get something like this:
firstly, if we fix a PER $R$, then $[FR \to R]$ is a PER of all algebras based on $R$. Lets say the limit of $U$ is $R_0$, then every element $f\in R_0$ restricts to a mapping $f_R:[FR \to R]\to R$. This is a morphism of PERs, because the category of PERs is a full subcategory of the effective topos. As a consequence $f_R$ itself is an element of the PER $[[FR \to R]\to R]$.
Now $[[FR \to R]\to R]_{R\in\PERs}$ is a family of PERs indexed by the class of PERs itself, and we can take $R_0\subset \prod_{R\in\PERs} [[FR \to R]\to R]$.

Secondly, the object of PERs exists within the effective topos, and is uniform. 
This makes $\bigcap_{R\in\PERs} [[FR \to R]\to R]$, the intersection of this family of PERs, already its product inside the category of PERs (see \cite{ms}). So to find $R_0$, we only need to select those elements of $\bigcap_{R\in\PERs} [[FR \to R]\to R]$ that commute with all the algebra morphisms. The results in the paper $\cite{din}$ seem to suggest that $R_0=\bigcap_{R\in\PERs} [[FR \to R]\to R]$. But in any case, the underlying PER of the initial algebra is:

\begin{equation}
R_0 := \left\{\begin{array}{r|l}
& \forall\ (R,a),(S,b),(T,c)\in F\textrm{-alg}_0,\\
(f,f')\in N^2 & \forall\ m:(R,a)\rightarrow (T,c), m':(S,b)\rightarrow (T,c).\\
& (m(fa) ,m'(f'b))\in T
\end{array}\right\}
\end{equation}

Given an algebra $(R,a)$, we get the projection map $\pi_a:R_0\to R$, simply defined by evaluation: $\pi_a(f) = fa$. These projections taken together form the limiting cone, which justifies calling $R_0$ the limit of $U$. Obviously, any algebra structure $c$ on $R_0$ has to make the following diagram commute for any algebra $(R,a)$:
$$
\xymatrix{ FR_0 \ar@{->}[r]^{F\pi_a} \ar[d]^c& FR\ar[d]^{a}\\
R_0 \ar[r]^{\pi_a} & R
}
$$
That means that for all $(x,y)\in FR_0$, $(cxa,a(F\pi_ay))\in R$. This can be achieved by letting $c:= \lambda xa.a(\phi\pi_ax)$ (for any $\phi$ tracking $F$).

The construction above shows that the category $\mathcal P$ of PERs is algebraically complete. The category $\mathcal C$ of pointed CEPERs is a reflective subcategory of $\mathcal P$, as is shown in $\cite{exp}$ (the prove of this fact doesn't use monotony, and is sound). Therefore it inherits weak completeness from $\mathcal P$, and we may use a similar construction to find an initial algebra for any endofunctor.
With the theory developed in \cite{rt} the fact that the category of pointed CEPERs is a CPO category can be used to prove that it is algebraically compact, \emph{because} it is algebraically complete.

\section{Yoneda}
Before we get to the Yoneda lemma, we need to know some things about natural transformations between realizable functors:

\begin{definition}\textnormal{
A natural transformation $\eta$ between two realizable endofunctors $F$ and $G$ of the category of PERs is realizable, iff there is a single number $e$ such that $\eta_R = [e]_{FR\to GR}$ for all PERs $R$.
}\end{definition}
Again, realizability is what makes the transformations internal to the effective topos. In this case the definition given in \cite{exp} is correct.

Because natural transformations are represented by natural numbers -- or because the category of PERs is weakly complete and internal: it al depends on your perspective -- we can construct a PER of natural transformations between any pair of PER valued functors. In fact: categories of PER valued functors are enriched over the category of PERs, as long as the domains are internal categories of the effective topos.

\begin{theorem}
Every endofunctor of $\PERs$ is naturally isomorphic to a monotone endofunctor.
\end{theorem}

\begin{proof} We know because of Yoneda's lemma that $FX\simeq \mathrm{nat}(\hom(X,-),F)$ naturally in both $F$ and $X$, and when $F$ is and endofunctor of $\PERs$, then the mapping $F_*$ satisfying $F_*X= \mathrm{nat}(\hom(X,-),F)$ can be turned into an endofunctor too: for any morphism $f:X\to Y$, if $\phi$ tracks $f$, then the induced morphism $f_*:\mathrm{nat}(\hom(X,-),F)\to \mathrm{nat}(\hom(X,-),F)$ is tracked by $\lambda xy.x\circ y\circ \phi$.
Now $F_*$ happens to be monotone:

If $X\subset Y$ and $[n]:Y\rightarrow Z$, then $[n]:X\rightarrow Z$ because $(nx,ny)\in Z$ whenever $(x,y)\in Y$ and $(x,y)\in Y$ whenever $(x,y)\in X$. Therefore $\hom(Y,-)\subset \hom(X,-)$ point wise. Furthermore: 
if $i$ is an index of the identity function, it determines a natural transformation $[i]:\hom(Y,-)\Rightarrow\hom(X,-)$.

Let $[i]:F\Rightarrow G$ for any two functors $F$ and $G$, and let $n:G\Rightarrow H$, then $n:F\Rightarrow G$, because $n\circ i$ represents the same p.r.f as $n$. Therefore $\mathrm{nat}(G,-)\subset \mathrm{nat}(F,-)$, and even $i:\mathrm{nat}(G,-)\Rightarrow \mathrm{nat}(F,-)$ point wise.

We see that if $X\subset Y$, then $F_*X\subset F_*Y$. Therefore $F$ is a monotone functor.
 \end{proof}

Although some internal functors may not be monotone, the assumption that realizable functors are, can be made without loss of generality. With this information added the original proof suffices to show that the category of pointed complete extensional PERs is indeed algebraically compact.

\end{document}